\documentclass{elsart}

\usepackage{amsfonts}

\usepackage{amsmath}

\usepackage{amssymb}
\usepackage{graphics}

\input{diagrams.tex}
\input{xy}
\xyoption{all}

\newcommand{\Abar}{\overline{A}}
\newcommand{\Am}{{A^{(m)}}}
\newcommand{\Ahat}{\hat{A}}
\newcommand{\Chat}{\hat{C}}
\newcommand{\Agrave}{{A^{^{\hspace{-2pt}\scriptscriptstyle\backprime}}}}
\newcommand{\Cgrave}{{C^{^{\hspace{-1pt}\scriptscriptstyle\backprime}}}}
\newcommand{\LAgrave}{\mathcal{L}_{{{A^{^{\hspace{-2pt}\scriptscriptstyle\backprime}}}}}}
\newcommand{\LmAgrave}{\mathcal{L}_m(A^{^{\hspace{-2pt}\scriptscriptstyle\backprime}})}
\newcommand{\Lm}{\mathcal{L}_m (\Am)}
\newcommand{\QA}{\mathcal{Q}_A}
\newcommand{\QC}{\mathcal{Q}_C}
\newcommand{\Hom}{\textup{Hom}\,}
\newcommand{\sHom}{\underline{\textup{Hom}}\,}
\newcommand{\Ext}{\textup{Ext}}

\newcommand{\ind}{\textup{ind}\,}
\newcommand{\add}{\textup{add}\,}

\newcommand{\Coker}{\textup{Coker}\,}
\newcommand{\pd}{\textup{pd}\,}
\newcommand{\isomorphe}{\cong}
\newcommand{\CA}{\mathcal{C}(A)}
\newcommand{\CmA}{\mathcal{C}_m(A)}

\newcommand{\LmC}{\mathcal{L}_m{{(C)}}}

\newcommand{\DA}{\mathcal{D}^b(\textup{mod}\, A)}
\newcommand{\DC}{\mathcal{D}^b(\textup{mod}\, C)}

\newcommand{\bimoduleC}{\,_C\!DC _C}

\newcommand{\ot}{\leftarrow}
\newcommand{\epito}{\twoheadrightarrow}
\newcommand{\monoto}{\hookrightarrow}

\newcommand{\zG}{\Gamma}
\newcommand{\zS}{\Sigma}
\newcommand{\zO}{\Omega}
\newcommand{\za}{\alpha}
\newcommand{\zb}{\beta}

\newcommand{\zl}{\lambda}

\begin{document}

\begin{frontmatter} 
  \title{$m$-cluster categories and $m$-replicated algebras} 
  \author{I. Assem\thanksref{nserc}}, \address{D\'epartement de
  Math\'ematiques, 
  Universit\'e de Sherbrooke, Sherbrooke (Qu\'ebec), J1K 2R1, Canada}
  \ead{ibrahim.assem@usherbrooke.ca}
  \author{T. Br\"ustle\thanksref{nserc+}}, \address{D\'epartement de
  Math\'ematiques, 
  Universit\'e de Sherbrooke, Sherbrooke (Qu\'ebec), J1K 2R1, Canada
  \emph{and}  
  Department of Mathematics, Bishop's University,  Lennoxville,
  (Qu\'ebec), 
  J1M 1Z7, Canada} 
  \ead{thomas.brustle@usherbrooke.ca}
  \author{R. Schiffler}, \address{Department of
  Mathematics and 
  Statistics, University of Massachusetts at Amherst, Amherst, MA
  01003-9305, USA}
  \ead{schiffler@math.umass.edu} 
  \author{G. Todorov}, \address{Department of Mathematics,
  Northeastern University, Boston, MA 02115, USA} 
  \ead{todorov.neu.edu}
  \thanks[nserc]{Partially supported by the NSERC of Canada}
  \thanks[nserc+]{Partially supported by the NSERC of Canada and the
  universities of Sherbrooke and Bishop's} 
\begin{abstract}
Let $A$ be a hereditary algebra over an algebraically closed field. 
We prove that an exact fundamental domain for the $m$-cluster category
$\CmA$ of 
$A$ is the $m$-left part $\mathcal{L}_m (A^{(m)})$ of the $m$-replicated
algebra of $A$.
Moreover, we obtain a one-to-one correspondence between the tilting
objects in $\CmA$ (that is, the $m$-clusters) and those tilting modules in
$\textup{mod}\, A^{(m)}$ for which all non projective-injective direct
summands lie in $\mathcal{L}_m (A^{(m)})$.

\noindent Furthermore, we study the module category of $A^{(m)}$ and
 show that a  basic exceptional module  with the correct number of
 non-isomorphic indecomposable summands is actually a tilting module. 
We also show how to
determine the projective  dimension of an indecomposable
$A^{(m)}$-module from its position in the Auslander-Reiten quiver.
 \end{abstract} 

\end{frontmatter} 

\setcounter{section}{-1}

\begin{section}{Introduction}
Cluster categories were introduced in \cite{BMRRT} and, for type
$A_n$, also in \cite{CCS} in order to understand better the cluster
algebras of Fomin and Zelevinsky \cite{FZ1,FZ2}. They are defined as
follows. Let $A$  be a hereditary algebra over an algebraically
closed field, and $\DA$ be the derived category of bounded complexes
of finitely generated $A$-modules, then the cluster category is the
orbit category of $\DA$ under the action of the functor $F=\tau^{-1} [1]$,
where $\tau$ is the Auslander-Reiten translation in $\DA$ and $[1]$ is
the shift. Later, the $m$-cluster category $\CmA$ was introduced in
\cite{T} (see also \cite{Keller-Reiten,Bin Zhu,Baur-Marsh}) 
 as a means for encoding  the combinatorics of $m$-clusters of
Fomin and Reading \cite{FR} in a fashion similar to the way the
cluster category encodes the combinatorics  of clusters.
     It is defined
to be the orbit category of $\DA$ under the  action of the functor
$\tau^{-1}[m]$. By \cite{K}, this category is triangulated. It is
proven in \cite{T} that there exists a bijection between $m$-clusters
and $m$-tilting sets in $\CmA$, that is, maximal sets of
indecomposables $S$ such that $\Ext^i_{\CmA}(X,Y)=0$ for all $X,Y$
in $S$ and all $i$ with $1\le i\le m$ (then the object $T=\oplus_{X\in
  S} X $ is called a {\em tilting object in } $\CmA$).

In \cite{ABST1},  we have given an interpretation of the cluster
category  and its tilting objects in terms of modules over a finite
dimensional algebra, namely the duplicated algebra of the original
hereditary algebra $A$. Our objective in the present paper is to
extend this characterization to the $m$-cluster category and its
tilting objects.

Following \cite{AI}, we define the $m$-replicated algebra of $A$ to be
the (finite dimensional) matrix algebra 
\begin{equation}\nonumber\Am=\left[
\begin{array}{ccccccccc}
\ A_0\ &\ 0 &\ldots& \ldots&\ldots&\ 0 \ \\ 
Q_1 & A_1 & 0 & \ldots &\ldots& 0\\  
0& Q_2 & A_2 & 0 & \ldots & 0\\
\vdots&&\ddots&\ddots&& \\
0&\ldots &&0& Q_m & A_m \\    
\end{array}\right]
\end{equation}
where $A_i=A,\ Q_i=DA$ for all $i$ and
all the remaining coefficients are zero (see \cite{HW} or section
\ref{section 1.4} below for the definition of the operations on
$\Am$).
Then $\Am $ is a quotient of the repetitive
algebra $\Ahat$ of $A$ (see \cite{HW}). Since $A$ is hereditary, the
structure of the module category $\textup{mod}\,\Am$ is known (see
section \ref{section 3.1} below). As a first useful consequence, we show that the
projective dimension of any indecomposable $\Am$-module is completely
determined by its position inside the module category (see Proposition
\ref{prop 3.4} below).

In order to relate the tilting $\Am$-modules to the tilting objects in
$\CmA$, we need to check whether exceptional modules with a maximal
number of summands are indeed tilting modules. We recall that, if $C$
is a finite dimensional algebra, a $C$-module $T$ is called
{\em exceptional} if 
\begin{enumerate}
\item the projective dimension $\pd T =d$ of $T$ is finite, and
\item $\Ext_C^i(T,T)=0 $ for all $i\ge 1$.
\\
An exceptional module $T$ is called a  {\em (generalized) tilting
  module} (see \cite{H}) if moreover :
\item there exists an exact sequence
\begin{equation*}
0\to C_C \to T_0 \to T_1 \to \cdots \to T_d \to 0
\end{equation*}
where each  $T_i$  is a direct sum of direct summands of $T$.
\end{enumerate}

It is an important open problem whether, for an exceptional module
$T$, 
having the number of isomorphism classes of indecomposable summands
equal to the rank of the Grothendieck group of $C$, is sufficient for
$T$ to be tilting.
This was first proven  by Bongartz in case $\pd T =1$, and the way he
did it was to prove that, if $T$ is exceptional, then there exists a
module $X$ such that $T\oplus X$ is a tilting module \cite{B}. This
latter statement (and hence the conjecture) were shown later for the
case where $C$ is representation-finite \cite{RS} (see also
\cite{CHU}). We prove here the analogue of Bongartz' result in another
special case and deduce our first theorem.

\begin{thm}\label{thmA}
Let $A$ be a hereditary algebra over an algebraically closed field
and $\Am$  be its     $m$-th replicated  algebra.
Let $T$ be a faithful exceptional $\Am$-module with $\pd T
\le m$, and such that the number of isomorphism classes of indecomposable
summands of $T$ equals the rank of the Grothendieck group. 
Then $T$ is a tilting module. 
\end{thm}

We then proceed to describe the $m$-cluster category $\CmA$. By
Happel's theorem \cite{H}, the derived category $\DA$ is equivalent to
the stable module category over the repetitive algebra $\Ahat $ of
$A$. The natural embedding of $\textup{mod}\,\Am$ into
$\textup{mod}\,\Ahat$ induces a functor $\pi$ from
$\textup{mod}\,\Am$ to the $m$-cluster category $\CmA$. Defining
the $m$-left part $\mathcal{L}_m(A^{(m)})$ of $\Am$ to consist of the
indecomposable 
$\Am$-modules all of whose predecessors have projective  dimension at most
$m$ and the functor $\pi$ to be the composition 
\begin{equation}\nonumber
 \pi: \textup{mod}\,\Am \monoto
\textup{mod}\,\Ahat\epito\underline{\textup{mod}}\,\Ahat \isomorphe
\DA \epito \CmA,
\end{equation}  
 this leads to the second theorem.

 \begin{thm} \label{thmB}
 Let $A$ be a hereditary algebra over an algebraically closed
 field   and $\Am$  be its     $m$-th replicated  algebra.
 The functor $\pi$ induces a one-to-one correspondence
 between the non projective-injective modules lying in the $m$-left part
 $\mathcal{L}_{m}(A^{(m)})$ and the indecomposable objects in 
 $\mathcal{C}_{m}(A)$. 
  \end{thm}

This is expressed by saying that $\mathcal{L}_{m}(A^{(m)})$ is an {\em
  exact fundamental domain}
for the functor $\pi$ (compare with \cite{ABST1}). We next
characterize the tilting objects in $\CmA$ in terms of the tilting
$\Am$-modules.  
 An exceptional $\Am$-module $T$ is called an
 $\mathcal{L}_{m}-exceptional$ module  if, whenever it is written in
 the form 
 $T=T'\oplus P$, with $P$ projective-injective, and $T'$ having no
 projective-injective direct summand, then all the indecomposable
 summands of $T'$ lie in $\mathcal{L}_m(A^{(m)})$.

  \begin{thm}\label{thmC}
   Let $A$ be a hereditary algebra over an algebraically closed field
   and $\Am$  be its $m$-th replicated  algebra.
   There is  a one-to-one correspondence between the
   $\mathcal{L}_{m}$-exceptional $\Am$-modules and the
   exceptional objects in $\CmA$ 
   given by $T=T'\oplus P \mapsto \pi (T').$ 
 \end{thm}
 
As a direct consequence of our Theorems \ref{thmA} and \ref{thmC}, the
above correspondence induces a one-to-one correspondence between the
$\mathcal{L}_m$-tilting $\Am$-modules and the tilting objects in $\CmA$.

Clearly, our theorems \ref{thmB} and \ref{thmC} generalize the main
results of \cite{ABST1}. The proofs here are however different, and
rest on the analysis of the projective  dimension of the modules under
consideration. 

We now describe the contents of the paper. After a brief preliminary
section, in which we fix the notations and recall the concepts needed
in the paper, our section 2 is devoted to the analysis of the
projective dimensions of the indecomposable modules over an infinite
dimensional (but locally finite dimensional) quotient of the
repetitive algebra, called the right repetitive algebra. These results
are then applied to the $m$-replicated algebra in section 3, which
culminates with the proof of Theorem \ref{thmB}. Section 4 is
devoted to the proof of Theorem \ref{thmA}, and section 5 to the proof
of Theorem \ref{thmC}.

\end{section}

\begin{section}{Preliminaries}
\begin{subsection}{Notation} 
Throughout this paper, algebras are basic and connected over a fixed
algebraically closed field. Given a locally finite dimensional
algebra $C$ (see \cite{BG}), we denote by $\textup{mod}\,C$ the category of finitely
generated right $C$-modules and by $\ind C$ a full subcategory whose
objects are a full set of representatives of the isomorphism classes of
indecomposable $C$-modules. Whenever we say that a given $C$-module is
indecomposable, we always mean implicitly that it belongs to $\ind
C$. Throughout this paper,  all subcategories of $\textup{mod}\,C$ are
full, and so are identified with their object classes. Given a
subcategory $\mathcal{C} $ of $\textup{mod}\,C$, we sometimes write
$M\in \mathcal{C}$ to express that $M$ is an object in
$\mathcal{C}$. We denote by $\add \mathcal{C}$ the subcategory of
$\textup{mod}\,C$ having as objects the finite direct sums of objects
in $\mathcal{C}$ and, if $M$ is a module, we abbreviate $\add\{M\}$ as
$\add M$. We denote the projective dimension of a module $M$ as $\pd
M$. The global dimension of $C$ is denoted by $\textup{gl.dim}\, C$ and
the quiver of $C$ by $\QC$.

Given an algebra $C$, we denote by $\nu_C= - \,\otimes_C DC$ its
Nakayama functor, and by $\tau_C$ its Auslander-Reiten
translation.
If $M$ is a $C$-module, then its first syzygy $\zO_{C}M$ is the
kernel of a projective cover $P\to M$ in $\textup{mod}\,C$ and its
first cosyzygy $\zO_{C}^{-1}M$ is the 
cokernel of an injective envelope $M\to I$.
 For further facts and definitions needed on
$\textup{mod}\,C$ and the Auslander-Reiten quiver $\zG(\textup{mod}\,C
)$ of $C$, we refer the reader to \cite{ARS,R}. For (minimal)
approximations we refer to \cite{AR}.
\end{subsection}

\begin{subsection}{The $m$-left part}
Let $C$ be a (locally) finite dimensional algebra, and $M,N$ be two
indecomposable $C$-modules. A {\em path} from $M$ to $N$ in $\ind C$ is
a sequence of non-zero morphisms
\begin{equation}\nonumber
  \label{eq*} M=M_0\stackrel{f_1\,}{\longrightarrow} M_1 \stackrel{f_2\,}{\longrightarrow}
  \cdots \stackrel{f_t\,}{\longrightarrow} M_t =N
\end{equation}
with all $M_i$ in $\ind C$. Following \cite{R}, we denote the
existence of such a path by $M\le N$. We
say that $M$ is a {\em predecessor} of $N$ (or that $N$ is a {\em successor}
of $M$). 

More generally,
if $S_1$ and $S_2$ are two sets of modules, we write $S_1\le S_2$
if every module in $S_2$ has a predecessor in $S_1$, every module in
$S_1$ has a successor in $S_2$ and no module in
$S_2$ has a successor in $S_1$  and no module in $S_1$ has a
predecessor in $S_2$. The notation $S_1 < S_2$ stands for 
$S_1\le S_2 $ and $S_1\cap S_2=\emptyset$. 

Let $m\ge 1$. 
The {\em $m$-left part} $\LmC$ of
$\textup{mod}\,C$ is the full subcategory of $\textup{ind}\,C$
consisting of all indecomposable $C$-modules $M$
such that if $L\le M$, then  $\pd L\le m$.

Clearly, $\mathcal{L}_1 (C)$ is the left part in the sense of \cite{HRS}.
\end{subsection}
  
\begin{subsection}{The cluster category and the $m$-cluster
    category}\label{section 1.3}
Let $A$ be a hereditary finite dimensional algebra, and $F$ denote the
endofunctor of $\DA$ defined as the composition $\tau^{-1}\,[1]$,
where $\tau $ is the Auslander-Reiten translation in $\DA$ and $[1]$
is the shift functor. The {\em cluster category } $\CA$ (see
\cite{BMRRT}) has as objects the $F$-orbits of objects in $\DA$ and
the morphisms are given by
\begin{equation*}
\Hom_{\CA}(\tilde X,\tilde Y)= \oplus _{i\in\Zset} \Hom_{\DA}(X,F^iY)
\end{equation*}
where $X$ and $Y$ are objects in $\DA$ and $\tilde X,\,\tilde Y$ are
their respective $F$-orbits. It is shown in \cite{K} that $\CA$ is a
triangulated category.

More generally, let $m\ge 1$ and $F_m$ denote the endofunctor of $\DA$
defined as the composition $\tau^{-1}\,[m]$. The $m$-{\em cluster
  category } $\CmA$ (see \cite{T}) has as objects the $F_m$-orbits of
objects in $\DA$ and the morphisms are given by
\begin{equation*}
\Hom_{\CmA}(\tilde X,\tilde Y)= \oplus _{i\in\Zset} \Hom_{\DA}(X,F_m^iY)
\end{equation*}
where $X$ and $Y$ are objects in $\DA$ and $\tilde X,\,\tilde Y$ are
their respective $F_m$-orbits. Again, by \cite{K}, $\CmA$ is a
triangulated category. We refer to \cite{T,Keller-Reiten,Bin Zhu,Baur-Marsh} 
for facts about the
$m$-cluster category.
\end{subsection}

\begin{subsection}{The repetitive algebra}\label{section 1.4}
Let $C$ be a finite dimensional algebra. Its 
repetitive algebra  $\Chat$ is  the infinite matrix algebra
\begin{equation}\nonumber\Chat  \quad = \quad \left[
  \begin{array}{cccccccccc}
  \ddots &&&\ 0\ \\
&\ C_{i-1}\ \\
&Q_i&\ C_i\ \\
&&Q_{i+1}&\ C_{i+1}\ \\
&\ 0\ &&&\ddots
  \end{array} \right]
\end{equation}
where matrices have only finitely many non-zero coefficients, $C_i=C$
and $Q_i= \bimoduleC$ for all $i\in \Zset$, all the remaining
coefficients are zero and multiplication is induced from the canonical
isomorphisms $C\otimes_C DC\isomorphe  \bimoduleC \isomorphe DC\otimes_C
C$ and the zero morphism $DC\otimes_C DC \to 0$, see \cite{HW}.
Then $\Chat $ is an infinite dimensional, locally finite dimensional,
self-injective algebra without identity. The importance of $\Chat$ in
our case comes from the following result of \cite{H}.
\begin{thm}\label{Happelthm}
\textup{(Happel)} Let $C$ be of
finite global dimension, then the derived category $\DC$ is
equivalent, as a triangulated category, to the stable module category
$\underline{\textup{mod}}\,\Chat$.
\end{thm}  
\end{subsection}

\end{section}

\begin{section}{The right repetitive algebra }

\begin{subsection}{Definition and description of the Auslander-Reiten
    quiver}\label{section 2.1} 
Let $C$ be a finite dimensional algebra. 
The {\em right repetitive algebra} $\Cgrave$ of $C$, introduced in
\cite{ABM}, 
is  the quotient of the repetitive algebra $\Chat $ of $C$ defined by:
\begin{equation}\nonumber
 \Cgrave \quad = \quad 
\left[ \begin{array}{cccccccccc}
\ C_0\ &\  && 0&& \\ 
Q_1 & C_1 &  &  && \\  
& Q_2 & C_2 &  &  & \\
0&&\ddots&\ddots&& \\
  \end{array}\right]
\end{equation}
where, as in section \ref{section 1.4},
$C_i=C$
and $Q_i= \bimoduleC$ for all $i\in \Zset$.

Assume, from now on, that $A$ is a hereditary algebra. The description
 of $\textup{mod}\,\Agrave$ follows easily from \cite{Ta,Y,HW,A} and
 can be summarized as follows.
\begin{lem}
\begin{enumerate}
\item The standard embeddings $\ind A_i\monoto \ind \Agrave$ (for
  $i\ge 0$) and $\ind \Agrave \monoto \ind \Ahat$ are full, exact,
  preserve indecomposable modules, almost split sequences and
  irreducible morphisms.
\item Under these embeddings, each $\ind A_i$ is a full convex
  subcategory of $\ind \Agrave$, and $\ind \Agrave$ is a full convex
  subcategory of $\ind \Ahat$. Moreover, $\ind A_0$ is closed under
  predecessors in $\ind \Agrave$, and $\ind \Agrave$ is closed under
  successors in $\ind \Ahat$.
\end{enumerate}
\end{lem}

In the sequel, we identify $A$ with $A_0$ and each $\ind A_i$ with the
 corresponding full subcategory of $\ind \Agrave$. Thus the
 Auslander-Reiten quiver $\zG(\textup{mod}\,\Agrave)$ of $\Agrave$
 can be described as follows (see Figure \ref{fig1}).

\begin{figure}
\includegraphics{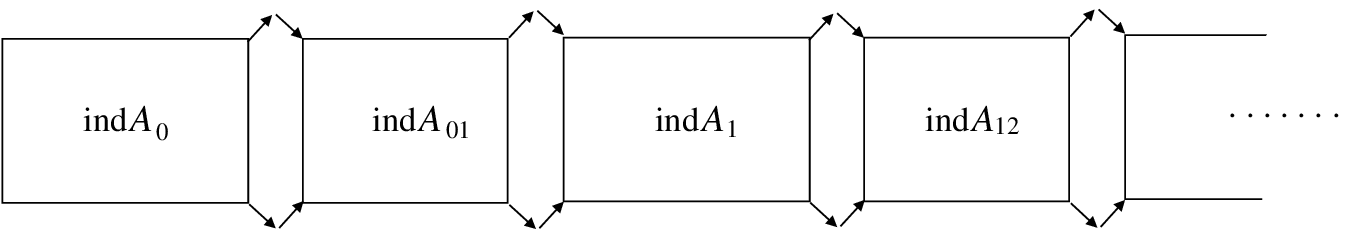}
\caption{Auslander-Reiten quiver of $\Agrave$}\label{fig1}
\end{figure}

It starts with  the Auslander-Reiten quiver  $\zG(\textup{mod}\,A)$
of $A_0=A$. Then pro\-jec\-tive-injective modules start to appear; such a
projective-injective module has its socle corresponding to a point in the
quiver $\QA$ of $A$, and its top corresponding to a point in the
quiver of $A_1$. Next occurs a part denoted by $\ind A_{01}$ where
indecomposables contain at the same time simple composition factors
from $A_0=A$, and simple composition factors
from $A_1$. When all projective-injectives whose socles correspond to
points in the quiver of $A_0$ have appeared, we reach the projective
$A_1$-modules and thus the Auslander-Reiten quiver
$\zG(\textup{mod}\,A_1)$ of $A_1$. The situation then repeats itself.

This repetition is effected by means of the Nakayama functor, whose
action is described as follows. For a point a (or an arrow $\za$) in the
ordinary quiver of $A$, denote by $a_i$ (or $\za_i$) the corresponding
point (or arrow, respectively) in the quiver of $A_i$. Let thus $M$ be
an $\Agrave$-module, considered as a representation. 
Then $\nu_{\Agrave}\, M(a_i)=M(a_{i-1})$ and
$\nu_{\Agrave}\, M(\za_i)=M(\za_{i-1}).$

Observe that all injective $\Agrave$-modules are projective, and that
those projective $\Agrave$-modules which are not injective are just
the projective $A$-modules.

\end{subsection}

\begin{subsection}{Example}

Since the former picture is the basis for our intuition, we give here
an example.

 Let $A$ be given by the quiver
\[\xymatrix{a\\&b\ar[lu]_\za\ar[ld]^\zb\\c}\]
then the quiver of $\Agrave $ is
\[
\xymatrix{a_0&&a_1\ar[dl]_{\zl_1}&&a_2\ar[dl]_{\zl_2}&&\cdots \\
&b_0\ar[ul]_{\za_0}\ar[dl]^{\zb_0}&&b_1\ar[ul]_{\za_1}\ar[dl]^{\zb_1}
&&b_2\ar[ul]_{\za_2}\ar[dl]^{\zb_2}&&\cdots\\
c_0&&c_1\ar[ul]^{\mu_1}&&c_2\ar[ul]^{\mu_2}&&\cdots
}
\]
bound by the relations $\zl_{i+1}\zb_i=0$, $\mu_{i+1}\za_i=0$,
$\za_{i+1}\zl_{i+1}=\zb_{i+1}\mu_{i+1}$, for all $i\ge 0$. Then
$\zG(\textup{mod}\,\Agrave)$ is as follows:

\[
\xymatrix@C=8pt@R=0pt{
&&&{\begin{array}{c}
\vspace{-10pt} c_1\\ \vspace{-10pt} b_0\\c_0 \end{array}}\ar[dr]
&&&&&&{\begin{array}{c} \vspace{-10pt} c_2\\ \vspace{-10pt} b_1\\ c_1 \end{array}}\ar[dr]\\
a_0\ar[dr]
&&{\begin{array}{c} \vspace{-10pt} b_0\\  c_0 \end{array}}\ar[ur]\ar[dr]
&&{\begin{array}{c} \vspace{-10pt} c_1\\  b_0 \end{array}}\ar[dr]
&&{\begin{array}{c} \vspace{-10pt} a_1 \end{array}}\ar[dr]
&&{\begin{array}{c} \vspace{-10pt} b_1\\  c_1 \end{array}}\ar[ur]\ar[dr]
&&{\begin{array}{c} \vspace{-10pt} c_2\\  b_1 \end{array}}\ar[dr]
&&{\begin{array}{c} \vspace{-10pt} a_2 \end{array}}
&\cdots
\\
&{\begin{array}{c} \vspace{-10pt} b_0\\  a_0\  c_0 \end{array}}\ar[ur]\ar[dr]
&&{\begin{array}{c}  b_0 \end{array}}\ar[ur]\ar[dr]
&&{\begin{array}{c} \vspace{-10pt}  a_1\  c_1\\ b_0 \end{array}}\ar[ur]\ar[dr]\ar[r]
&{\begin{array}{c} \vspace{-10pt} b_1\\ \vspace{-10pt} a_1\  c_1\\ b_0 \end{array}}\ar[r]
&{\begin{array}{c} \vspace{-10pt} b_1\\  a_1\  c_1 \end{array}}\ar[ur]\ar[dr]
&&{\begin{array}{c}  b_1 \end{array}}\ar[ur]\ar[dr]
&&{\begin{array}{c} \vspace{-10pt}  a_2\  c_2\\ b_1 \end{array}}\ar[ur]\ar[dr]\ar[r]
&{\begin{array}{c} \vspace{-10pt} b_2\\ \vspace{-10pt} a_2\  c_2\\ b_1 \end{array}}
&\cdots
\\
c_0\ar[ur]
&&{\begin{array}{c} \vspace{-10pt} b_0\\   a_0 \end{array}}\ar[ur]\ar[dr]
&&{\begin{array}{c} \vspace{-10pt} a_1\\   b_0 \end{array}}\ar[ur]
&&{\begin{array}{c}  c_1\\    \end{array}}\ar[ur]
&&{\begin{array}{c} \vspace{-10pt} b_1\\   a_1 \end{array}}\ar[ur]\ar[dr]
&&{\begin{array}{c} \vspace{-10pt} a_2\\   b_1 \end{array}}\ar[ur]
&&{\begin{array}{c}  c_2\\   \end{array}}
&\cdots
\\
&&&{\begin{array}{c} \vspace{-10pt} a_1\\ \vspace{-10pt} b_0\\ a_0 \end{array}}\ar[ur]
&&&&&&{\begin{array}{c} \vspace{-10pt} a_2\\ \vspace{-10pt} b_1\\ a_1 \end{array}}\ar[ur]\\
}
\]
where modules are represented by their Loewy series.

\end{subsection}

\begin{subsection}{Injective envelopes and projective covers} 
The following lemma is inspired from a well-known
result about symmetric algebras.

\begin{lem}\label{lemma1}
Let $L$ be an indecomposable  $\Agrave$-module and 
\begin{equation*}
0\to L \stackrel{f}{\to} I \stackrel{g}{\to} N\to 0  
\end{equation*} 
 a short
exact sequence of $\Agrave$-modules
with $L \stackrel{f}{\to} I$ an injective envelope, $I$
projective-injective and $N\ne 0$.
Then 
\begin{enumerate}
\item $I\stackrel{g}{\to}N$ is a projective
cover in $\textup{mod}\,\Agrave$ and 
\item $N$ is indecomposable. 
\end{enumerate}
\end{lem}
\begin{pf} (1)
  Let $p:P(N)\to N$ be a projective cover in
  $\textup{mod}\,\Agrave$ and  let $K=\textup{Ker}\, p$. Then there exists a
  projective $\Agrave$-module $P' $ and a commutative diagram
  with exact rows and columns, where $i$ is a section.
  \begin{diagram}[size=2em]
  &&0&&0&&\\
  &&\dTo&&\dTo&&\\
  0&\rTo& K& \rTo &P(N)&\rTo^p & N& \rTo &0 &&\\
  &&\dTo&&\dTo^i && \dTo^{\textup {id}}\\
  0&\rTo& L& \rTo^f &I&\rTo^g& N& \rTo &0 &&\\
  &&\dTo&&\dTo\\
  &&P'&\rTo^{\textup {id}}&P'\\
 &&\dTo&&\dTo\\
  &&0&&0\\
\end{diagram}
Hence $P'$ is also injective.
Then $L\isomorphe K\oplus P'$.
Since $L$ is indecomposable and 
 not projective-injective it follows that
$L\isomorphe K$. Thus $P(N)=I$. 

(2) Suppose that 
$N=N_1\oplus N_2$. Then $I\isomorphe P(N)\isomorphe P(N_1)\oplus P(N_2)$
and  there exists an induced direct sum decomposition of the kernel
$L=L_1\oplus L_2 $. 
Since $L$ is indecomposable, we may assume that $L=L_1$ and $L_2=0$.
But then $P(N_2)\isomorphe
N_2$
so that we have a commutative diagram  with exact rows. 
\begin{diagram}[size=2em]
0&\rTo& L &\rTo^f& I& \rTo^g& N &\rTo& 0\\
&&\dTo^{\textup {id}}&&\dTo^\isomorphe&&\dTo^{\textup {id}}\\
0&\rTo& L &\rTo& P(N_1)\oplus N_2&
\rTo^{\left({{\scriptstyle k}\ {\scriptstyle 0}\atop 
 {\scriptstyle  0}\  {\scriptstyle \rm id}}\right)}& N_1\oplus N_2
&\rTo&  0\\ 
\end{diagram}
 which contradicts  the hypothesis  that
 $L\stackrel{f}{\to}I$ is an injective envelope.  
\qed
\end{pf}

\begin{cor}\label{lemma1'}
Let $L$ be an indecomposable  $\Agrave$-module and
\[\xymatrix{
0\ar[r]&
L\ar[r]^{f^0} &
I^0\ar[r]^{f^1} &
I^1\ar[r]&
\ldots\ar[r]^{f^k} &
I^k
}
\]
%
be a minimal injective coresolution in $\textup{mod}\,\Agrave$, with all $I^j$
projective-injective and $N=\Coker f^k\ne 0$.
Then:
\begin{enumerate}
\item $\xymatrix{
I^0\ar[r]^{f^1} &
I^1\ar[r]&
\ldots\ar[r]^{f^k} &
I^k\ar[r]&
N\ar[r]&
0
}$
 is a minimal projective
resolution in $\textup{mod}\,\Agrave$,
\item $N$ is indecomposable and 
\item for all $j \le k$,
the $j$-th cosyzygy of $L$ is isomorphic to the
$(k+1-j)$-th syzygy of $N$, that is 
\begin{equation}\nonumber
\zO^{-j}_{\Agrave}\,L\  \isomorphe \  \zO^{(k+1)-j}_{\Agrave}\,N.
\end{equation}
\end{enumerate}
\end{cor}
\begin{pf}
  By induction on $k$ using Lemma \ref{lemma1}.\qed
\end{pf}

\begin{lem}\label{lem5}
Let $M$ be an indecomposable $\Agrave$-module which does not lie in $\ind
A$. Then the projective cover of $M$ is projective-injective.
\end{lem}

\begin{pf} We may clearly assume that
 $M$ is not projective-injective. 
Since $M\notin \ind A$, we have
\begin{equation}\nonumber
\overline\Hom_{\Agrave}(A_A,M)\isomorphe\overline\Hom_{\Ahat}(A_A,M)=0 ,
\end{equation}
 where the last equality follows from
Happel's Theorem (see  Theorem \ref{Happelthm} above) and from the structure of morphisms in the
derived category. Therefore, any non-zero morphism in
$\Hom_{\Agrave}(A_A,M)$ must factor through an injective $\Agrave$-module
which is also projective. The statement follows. \qed
\end{pf}

\begin{lem}
Let $M$ be an indecomposable $\Agrave$-module, then 
\begin{equation}\nonumber
\zO_{\Agrave}^{-1}\,\tau_{\Agrave}^{-1}\,M\ \isomorphe\
\tau_{\Agrave}^{-1}\,\zO_{\Agrave}^{-1}\, M
\end{equation}
\end{lem}
\begin{pf}
For any indecomposable $\Agrave$-module $M$, we have
$\zO_{\Agrave}^{-1}\,M=\zO_{\Ahat}^{-1}\,M$, because 
 injective $\Agrave$-modules are also $\Ahat$-injective,
 and
$\tau_{\Agrave}^{-1}\,M = \tau_{\Ahat}^{-1}\,M$. The statement follows
from the fact that $\tau_{\Ahat}^{-1}=\zO_{\Ahat}^{-2}
\nu_{\Ahat}^{-1}$ (see \cite[IV.3.7 p.126]{ARS}).
\qed
\end{pf}
\end{subsection}

\begin{subsection}{Projective dimension }
We are now able to prove the main result of this section.

\begin{thm}\label{thm1}
Let $M$ be an indecomposable $ \Agrave$-module and let $k\ge 1$. Then
$\pd M = k$  if and
    only if there exists an indecomposable $A$-module $N$ such that 
$ M \isomorphe \tau^{-1}_{\Agrave} \,\zO^{-(k-1)}_{\Agrave} \,N$.
\end{thm}

\begin{pf} We prove the statement 
by induction on $k$.
Suppose $k=1$.  If
$M \isomorphe \tau_{\Agrave}^{-1}\,N$ for some $N$ in $\ind A$
 then it follows from \cite[Cor 5 and  6]{ABST1} that 
 $M $ lies in the left part $ \LAgrave$ of $\Agrave$, and thus, $\pd M \le 1$.
On the other hand, $M\isomorphe\tau_{\Agrave}^{-1}\,N$ implies that $M$ is not
projective, and hence $\pd M =1$.

Conversely, assume that $M$ is not of the form
$\tau_{\Agrave}^{-1}\,N$, for some $N\in \ind A$. If $M$ is
projective then $\pd M=0\ne 1$. Otherwise, by Lemma \ref{lem5}, in 
the short exact sequence 
\begin{equation}\label{eqthm5}
0\to\zO_{\Agrave}\,M\to P\stackrel{f}{\to} M \to 0,
\end{equation}
where $f$ is a projective cover, the module $P$  is
projective-injective and, by \cite[Prop. 7]{ABST1},
$\zO_{\Agrave}\,M$ is not a projective $A$-module. 
Clearly, $\zO_{\Agrave}\,M$ cannot be projective-injective,
 since (\ref{eqthm5}) is not split. Hence $\pd M\ne 1$.

 Suppose now $k>1$.
Let $M\in \ind \Agrave$ such that $\pd M=k $. 
Since $k>1$, then $M\notin\ind  A$ and, by Lemma \ref{lem5}, the projective
cover of $M$ in $\textup{mod}\,\Agrave$ is
projective-injective. 
Let $X=\zO_{\Agrave}\, M $. Then 
$\pd X=k-1$. By the dual of Lemma \ref{lemma1},
 $X$ is indecomposable.
Now, by the induction hypothesis, $X= \tau^{-1}_{\Agrave}\,Y$ for some
 $Y\in\zO^{-(k-2)}_{\Agrave}\,\ind A$. 
Moreover, by Lemma \ref{lemma1}, $M=\zO^{-1}_{\Agrave}\,X$  and therefore 
\begin{equation}\nonumber
M=\zO^{-1}_{\Agrave}\,X=\zO^{-1}_{\Agrave}\,\tau^{-1}_{\Agrave}\,Y \isomorphe
\tau^{-1}_{\Agrave}\,\zO^{-1}_{\Agrave}\,Y
\end{equation} 
lies in $ \tau^{-1}_{\Agrave}\,\zO^{-(k-1)}_{\Agrave}\,\ind A$, as required.

Conversely, suppose that
$M=\tau^{-1}_{\Agrave}\,\zO^{-(k-1)}_{\Agrave}\,N$ for some $N\in \ind
A$. By Corollary \ref{lemma1'},
$X=\tau^{-1}_{\Agrave}\,\zO^{-(k-2)}_{\Agrave}\,N$
 is indecomposable.  
Furthermore, $M=\zO^{-1}_{\Agrave}\,X$
and there is a short exact sequence 
\begin{equation}\nonumber
0\to X
\stackrel{f}{\to} I \stackrel{g}{\to} M \to 0,
\end{equation} 
with $f$  an injective envelope in $\textup{mod}\,\Agrave$.
 Thus $I$ is projective-injective  and, by Lemma \ref{lemma1}, $g$ is
 a projective cover; hence $\pd M= \pd X +1$. 
But by induction, the projective dimension of
$X$ is  $k-1$, and thus  $\pd M=k$.
This completes the proof.\qed
\end{pf}

Recall that a component of the Auslander-Reiten
quiver of $\Agrave$ is called {\em trans\-jective} if it does not
contain cycles in $\ind \Agrave$.
A full connected subquiver $\zS$ in an Auslander-Reiten component is
called a subsection if every path in $\zS$ is sectional.
\begin{defn}
A \emph{right stable slice} $\zS$ in $\ind\Agrave$ is a connected
convex subsection in  a transjective component of the Auslander-Reiten
quiver of $\textup{mod}\,\Agrave$ which intersects each right stable
$\tau_{\Agrave}$-orbit in that component.
\end{defn}

For each $k$ greater than or equal to zero, we
 set  $\zS_k = \left\{
  \zO_{\Agrave}^{-k}\,P_x\mid x\in(\QA)_0\right\}$.
 Notice that $\zS_0$ is just the set of projective
$ A$-modules and that $\zS_k = \zO^{-k} (\zS_0)$.

\begin{lem}\label{slice}
 For each $k\ge 0$,  $\zS_k$ is 
  a right stable slice in $\ind\Agrave$.
\end{lem}
\begin{pf} By induction on $k$. Clearly
$\zS_0$,which is the set of all indecomposable projective $A$-modules,
  is a right 
  stable slice. Also, by \cite[Cor. 6]{ABST1},
  $\zS_1=\{\zO_{\Agrave}^{-1}\,P_x\mid x\in (\QA)_0\}$ 
  is a right stable slice.
Assume now that $k>1$. If $k$ is even we have $ \zO_{\Agrave}^{-k}=
  \zO_{\Agrave}^{-2(k/2)}=
  \tau_{\Agrave}^{-k/2}\,\nu_{\Agrave}^{k/2}$, by \cite[IV.3.7 p.126]{ARS}.
Since both
  $\tau_{\Agrave}$ and $\nu_{\Agrave}$ preserve right stable slices,
  the statement follows from
  $\zS_k=\tau_{\Agrave}^{-k/2}\,\nu_{\Agrave}^{k/2}\,\zS_0$. If, on
  the other hand,  $k $ is
  odd, then the statement follows similarly from 
  $\zS_k=\tau_{\Agrave}^{-(k-1)/2}\,\nu_{\Agrave}^{(k-1)/2}\,\zS_1$.
 \qed
\end{pf}

We now show that the right stable slices $\zS_k$ 
 partition  $\underline{\textup{ind}}\,\Agrave$ into regions of
constant projective dimension. It is useful to note that, if $L\le M$
in $\ind\Agrave$, then $\zO_{\Agrave}\,L \le\zO_{\Agrave}\, M$ and  
$\zO^{-1}_{\Agrave}\,L \le\zO^{-1}_{\Agrave}\, M$. 
\begin{cor}\label{corpd}
Let $M$ be an indecomposable $\Agrave$-module which is not
projective-injective and let $k\ge 1$. Then
\begin{equation}\nonumber
 \pd M =k \quad \text{ if and only if } \quad \zS_{k-1} < M \le \zS_k
\end{equation}
\end{cor}
\begin{pf} Suppose $\pd M =k$. By Theorem \ref{thm1}, there exists
  $N\in \ind A$ such that 
\begin{equation}\nonumber
M=\tau^{-1}_{\Agrave}\,\zO^{-(k-1)}_{\Agrave}\,N =
    \zO^{-(k-1)}_{\Agrave}\,\left(\tau^{-1}_{\Agrave}\,N\right) .
\end{equation}  
Now, by \cite{ABST1}, we have $\zS_0<\tau^{-1}_{\Agrave}\,N\le \zS_1$.
    Since  $\zS_{k-1}=\zO^{-(k-1)}_{\Agrave}\,\zS_0$ and
    $\zS_{k}=\zO^{-(k-1)}_{\Agrave}\,\zS_1$, by
 definition,  this implies $ \zS_{k-1}
    <M\le\zS_{k}$. 

Conversely, assume that $ \zS_{k-1} <M\le\zS_{k}$. Since, by Corollary
\ref{lemma1'}, we have $\zS_0 =\zO^{k-1}_{\Agrave} \zS_{k-1}$ and
$\zS_1 =\zO^{k-1}_{\Agrave} \zS_{k}$, applying
$\zO^{k-1}_{\Agrave}$ to the former inequality yields 
$\zS_0<\zO^{k-1}_{\Agrave}\,M\le \zS_1$.
 Then, by
\cite[Cor. 5]{ABST1}, it follows that $\pd \zO^{k-1}_{\Agrave}\,M =1$,
whence $\pd M=k$.\qed
\end{pf}

The right repetitive algebra $\Agrave$ thus provides an example of an
(infinite-dimensional) algebra such that every indecomposable
$\Agrave$-module has finite projective dimension, but such that the
finitistic projective dimension of $\Agrave $ is infinite.
  (We recall that
the finitistic projective dimension of an algebra is the supremum of the
projective dimensions of those modules which have finite projective
dimension.)

\end{subsection}
\end{section}

\begin{section}{The $m$-replicated algebra $ \Am $} 
\begin{subsection}{Definition and description of the Auslander-Reiten quiver}\label{section 3.1}
Let $m\ge 1$ and $C$ be a finite dimensional algebra. 
The \emph{$m$-replicated algebra} $C^{(m)}$  of  $C$ is  the quotient
of the right repetitive algebra $\Cgrave$ (hence of the repetitive
algebra $\Chat$) defined by
 \begin{equation}\nonumber C^{(m)}=\left[
\begin{array}{ccccccccc}
\ C_0\ &\ 0 &\ldots& \ldots&\ldots&\ 0 \ \\ 
Q_1 & C_1 & 0 & \ldots &\ldots& 0\\  
0& Q_2 & C_2 & 0 & \ldots & 0\\
\vdots&&\ddots&\ddots&& \\
0&\ldots &&0& Q_m & C_m \\    
\end{array}\right].
\end{equation}
If $m=1$, then $C^{(1)}$ is called the {\em duplicated algebra} of $C$ (see
\cite{ABST1}). It is shown in \cite{AI} that
\begin{equation}\nonumber
m+\textup{gl.dim}\,C\le\textup{gl.dim}\,C^{(m)}\le
(m+1)\textup{gl.dim}\,C+m.
\end{equation}
Assume from now on that $A$ is hereditary. The description of $\Am$
follows from \cite{Ta,Y,A} and can be summarized as follows.

\begin{lem}\label{lem 3.1}
\begin{enumerate}
\item The standard embeddings $\ind A_i\monoto \ind \Am$ (where $0\le
  i\le m $) and $\ind \Am \monoto \ind \Agrave$ are full, exact,
  preserve indecomposable modules, almost split sequences and
  irreducible morphisms.
\item Under the above embeddings, each $\ind A_i$ is a full convex
  subcategory of $\ind \Am$, and $\ind \Am$ is a full convex
  subcategory of $\ind \Agrave$. Moreover, $\ind A_0$ is closed under
  predecessors  and $\ind A_m$ is closed under
  successors in $\ind \Am$.
\end{enumerate}
\end{lem}

In the sequel, we identify $A$ with $A_0$ and each $\ind A_i$ with the
 corresponding full subcategory of $\ind \Am$, and $\ind \Am$ with the
 corresponding full subcategory of $\ind \Agrave$. The
 Auslander-Reiten quiver  of $\Am$
 can be deduced from that of $\Agrave$ (see Figure \ref{fig2}).

\begin{figure}
\scalebox{0.99}{\includegraphics{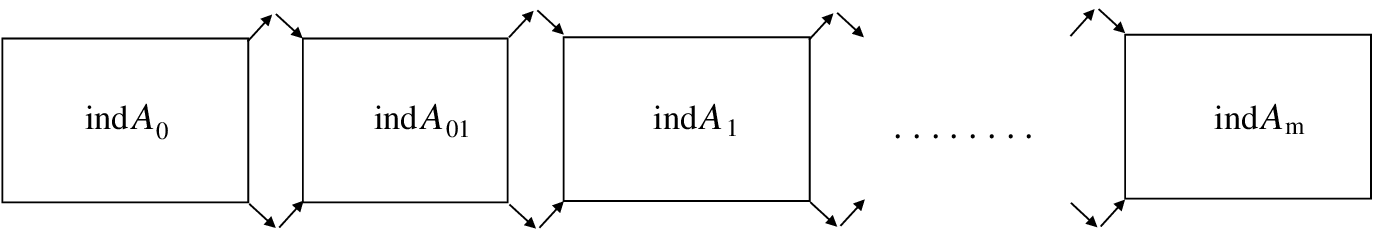}}
\caption{Auslander-Reiten quiver of $\Am$}\label{fig2}
\end{figure}

\end{subsection}

\begin{subsection}{Projective dimension }
Since we are interested in the projective dimension of indecomposable
$\Am$-modules, we need to study projective covers. In our situation,
Lemma \ref{lem5} becomes the following lemma.
\begin{lem}
\begin{enumerate} 
\item Let $M$ be an indecomposable $\Am$-module which does not lie in $\ind
A$. Then its projective cover in $\textup{mod}\,\Am$ is
projective-injective and coincides with its projective cover in
$\textup{mod}\,\Agrave$. 
\item Let $M$ be an indecomposable $\Am$-module which does not lie in
$\ind A_m$. Then its injective envelope in $\textup{mod}\,\Am$ is
  projective-injective and coincides with its injective envelope in
$\textup{mod}\,\Agrave$.
\end{enumerate}
\end{lem}

\begin{pf} Since (2) is dual to (1),  we only
 prove (1).
The case where $M$ is projective-injective is trivial. Assume that
$M$ is injective but not projective. Then $M\in\ind A_m$ and therefore
$\Hom_{\Am}(A,M)=0$ (because the supports of $A$ and $M$ are disjoint)
and the statement follows. 
In the case where  $M$ is not injective one can use the same
argument as in the proof of Lemma \ref{lem5} to show that the
projective cover of $M$ is projective-injective. It coincides with its
projective cover in $\textup{mod}\,\Agrave$ because $\ind \Am$ is closed under
predecessors in $\ind \Agrave$.
\qed
\end{pf}

We now show that, for each $k \le m$ the right stable slice $\zS_k$ of
Lemma \ref{slice} consists of $\Am$-modules. 
\begin{lem}\label{lemma 3.3}
We have $\zS_m\subset \ind \Am$.
\end{lem}
\begin{pf}
We prove by induction on $k\le m$ that $\zS_k\le\nu_{\Agrave}^k
\zS_0$. If $k= 1$, then by \cite{ABST1},
$\zS_1=\zO_{\Agrave}^{-1}\,\zS_0 \le \nu_{\Agrave}\,\zS_0$. Assume now
that, for some $k<m$, we have $\zS_k\le \nu_\Agrave^k\,\zS_0$. Since 
$\nu_{\Agrave} = \tau_{\Agrave}\,\zO_{\Agrave}^{-2}$, we have 
\begin{equation}\nonumber
 \nu_{\Agrave}^{k+1}\,\zS_0\ = \ \tau_{\Agrave}\,\zO_{\Agrave}^{-2}\,
 \nu_{\Agrave}^{k}\,\zS_0 \ \ge\  \tau_{\Agrave}\,\zO_{\Agrave}^{-2}
 \,\zS_k 
\ =\  \tau_{\Agrave}\,\zO_{\Agrave}^{-1} \,\zS_{k+1}.
\end{equation}
On the other hand, if $X\in\ind\Agrave$ is not a projective-injective module
then $\Hom_{\Agrave}(\tau_{\Agrave}^{-1}\,X,\zO_{\Agrave}^{-1}\,X)\ne
0$ so that $\tau_{\Agrave}^{-1}\,X \,\le \,\zO_{\Agrave}^{-1}\,X$ and
then $X\le \tau_{\Agrave}\,\zO_{\Agrave}^{-1}\,X$. Therefore
$\nu_{\Agrave}^{k+1}\,\zS_0\ge \zS_{k+1}$, establishing our claim. The
statement of the lemma then follows from
$\zS_m\le\nu_{\Agrave}^m\,\zS_0={A_m}$.
\qed
\end{pf}

We are now able to prove the statement
  corresponding to Theorem \ref{thm1} and Corollary \ref{corpd} in
  $\textup{mod}\,\Am$.

\begin{prop} \label{prop 3.4}Let $M$ be an indecomposable  $\Am$-module which is not projective and
let $k$ be such that $1\le k \le m$. The following are equivalent:
\begin{enumerate}
\item $\pd M =k ,$
\item there exists $N\in \ind A$ such that $M\isomorphe
  \tau_\Am^{-1}\,\zO_\Am^{-(k-1)}\,N,$
\item $\zS_{k-1}<M\le\zS_k$.
\end{enumerate}  
\end{prop} 

\begin{pf}
Since $\ind \Am$ is closed under predecessors in $\ind \Agrave$, then a
minimal projective resolution of an indecomposable $\Am$-module $X$ in
$\textup{mod}\,\Am$ is also a minimal projective resolution of $X$ in
$\textup{mod}\, \Agrave$ (and, in particular, the respective syzygies
coincide). Therefore, by Theorem \ref{thm1}, (1) holds if and only if
there exists $N\in \ind A $ such that  
\begin{equation}\nonumber
M\isomorphe \tau_{\Agrave}^{-1}\,\zO_{\Agrave}^{-(k-1)}\,N
\isomorphe  \tau_{\Am}^{-1}\,\zO_{\Am}^{-(k-1)}\,N.
\end{equation}  
Finally, the equivalence of (1) and (3) follows from Lemma \ref{lemma
  3.3} and Corollary \ref{corpd}.
\qed
\end{pf}

\begin{cor}\label{cor 3.5}
Let $M$ be an indecomposable $\Am$-module which is not
pro\-jec\-tive-injective. The following are equivalent:
\begin{enumerate}
\item $M\in \Lm$,
\item $M \in \cup_{k=1}^m \zO_\Am^{-(k-1)}(\ind A)\,\cup\,\zS_m$,
\item $M\le \zS_m$.
\end{enumerate}   
\end{cor}   

\begin{pf}
The equivalence of (1) and (3) follows from Proposition \ref{prop
  3.4}. We show that (1) implies (2). Let $M\in \Lm$, then $\pd
M=k\le m$. By Proposition \ref{prop 3.4}, there exists $N\in \ind A$
such that 
\begin{equation}\nonumber
M\,=
\, \tau_{\Agrave}^{-1}\,\zO_{\Agrave}^{-(k-1)}\,N
\,\isomorphe\, 
\zO_{\Agrave}^{-(k-1)}\,\tau_{\Agrave}^{-1}\,N
\,\isomorphe\, 
\zO_{\Am}^{-(k-1)}\,\tau_{\Am}^{-1}\,N.
\end{equation}  
Thus, if $N$ is not an injective $A$-module, then $M\in
\zO_\Am^{-(k-1)}(\ind A)$. On the other hand, if $N$ is an injective
$A$-module, then $\tau_\Am^{-1}\,N\in \zO_\Am^{-1}\,\zS_0$ (by
\cite{ABST1}) so that $M\in
\zO_\Am^{-(k-1)}\,\zO_\Am^{-1}\,\zS_0\,=\,\zO_\Am^{-k}\,\zS_0\,=\,\zS_k$.
Thus $M\in\zS_m$ if $k=m$ and $M\in \cup_{k=1}^m \zO_\Am^{-(k-1)}(\ind
A)$ otherwise.

Finally we prove that (2) implies (3). We assume $M\in\cup_{k=1}^m
\zO_\Am^{-(k-1)}(\ind 
A)\,\cup\,\zS_m$ and claim that $M\le \zS_m$. If $M\in \zS_m$, there
is nothing to prove. Otherwise, there exist $N\in \ind A$ and $k\le
m$ such that $M=\zO_\Am^{-(k-1)}\,N$. Now $N\le \zS_1 = \zO_\Am^{-1}\zS_0$
so $M\le
\zO_\Am^{-(k-1)}\,\zO_\Am^{-1}\,\zS_0\,=\,\zO_\Am^{-k}\,\zS_0\,=\,\zS_k \le
\zS_m$. 
\qed
\end{pf}

\end{subsection}

\begin{subsection}{Exact fundamental domain}
One very easy consequence of Corollary \ref{cor 3.5}, Theorem
\ref{thm1} and Corollary \ref{corpd} is that $\Lm=\LmAgrave$. Another
one is the following corollary.

\begin{cor}\label{cor 3.6}
The embedding functor 
\[\add\Lm \monoto \textup{mod}\,\Am \monoto
\textup{mod}\, \Agrave \monoto
\textup{mod}\,\Ahat \]
  is full, 
  exact and preserves indecomposable modules, irreducible morphisms
  and almost split sequences.
\end{cor}

We are now able to prove Theorem \ref{thmB}. Since $A$ is hereditary, we have,
by \cite{H}, an equivalence
$\DA\isomorphe\underline{\textup{mod}}\,\Ahat$ of triangulated
categories. Let
\begin{equation}\nonumber
 \pi: 
\textup{mod}\,\Am \monoto
\textup{mod}\,\Agrave \monoto
\textup{mod}\,\Ahat\epito\underline{\textup{mod}}\,\Ahat \isomorphe
\DA \epito \CmA.
\end{equation}  
be the canonical functor (where $\CmA$ is the $m$-cluster category,
see section \ref{section 1.3}).
We define an {\em exact fundamental domain} for
  $\pi$ to be a full convex subcategory of $\ind
  \Ahat$ which contains exactly one point of each fibre $\pi
  ^{-1}(X)$, with $X$ an indecomposable object in $ \CmA$.

\begin{thm}\label{thm2} 
The functor $\pi$ induces a one-to-one correspondence between
the non projective-injective modules in $\Lm$ and the indecomposable
objects in $\CmA$. In particular, $\Lm$ is an exact fundamental domain
for $\pi$. 
\end{thm}
\begin{pf} 
 Since
$\Lm$ is a full convex subcategory of $\ind \Am$ and $\ind\Agrave$, it
is also convex inside $\ind \Ahat$. Furthermore, the non
projective-injective modules in $\Lm$ are just the modules in 
\begin{equation}\nonumber
\cup_{k=1}^m \zO^{-(k-1)}_{\Am}\,(\ind A)\cup \zS_m
\ =\ 
\cup_{k=1}^m \zO^{-(k-1)}_{\Ahat}\,(\ind A)\cup \zS_m.
\end{equation}
 The
statement follows from the definition of $\CmA$ and the fact that
$\zO^{-1}_{\Ahat}$ corresponds to the shift $[1]$ under the equivalence
$ \underline{\textup{mod}}\,\Ahat \isomorphe\DA$.\qed
\end{pf}

\end{subsection}

\end{section}

\begin{section}{Tilting modules over the replicated algebra}

\begin{subsection}{Definitions and preparatory results}
Let $C$ be a finite dimensional algebra, and $T$ be a $C$-module. We
say that $T$ is an {\em exceptional module} if
\begin{enumerate}
\item $\pd T =d<\infty$
\item $\Ext_C^i(T,T)=0$ for all $i\ge 1$\\
An exceptional module $T$ is called a  {\em (generalized)  tilting
  module} if moreover
\item there exists an exact sequence
\begin{equation*}
0\to C_C \to T_0 \to T_1 \to \cdots \to T_d \to 0
\end{equation*}
where each  $T_i\in \add T$ for all $i$.  
\end{enumerate}   

It is useful to observe that, if $T$ is an exceptional module which is
faithful, then any projective-injective indecomposable $C$-module $P$
is a direct summand of $T$: indeed, since $T$ is faithful, there
exists a monomorphism $C_C\monoto T_0$ with $T_0 \in \add T$, which,
when composed with the inclusion $P\monoto C_C$ yields an inclusion
$P\monoto T_0$ which splits, because $P$ is injective.
In particular, if $ T$ is a tilting $C$-module, then any
projective-injective indecomposable  $C$-module is a direct summand of
$T$.

An exceptional module $T$ is said to be {\em basic} (or {\em
  multiplicity free}) if, whenever $T=\oplus_{i=1}^n T_i$ where all the
  $T_i$ are indecomposable, we have $T_i\ne T_j$ for $i\ne j$. It is
  well-known that, if $T$ is a basic tilting $C$-module, then the
  number of its indecomposable summands is equal to the rank of the
  Grothendieck group of $C$, see \cite{H}.

For the rest of this section, we let as before $A$ be a hereditary
algebra and $m\ge 1$, and $\Am$ denote the $m$-replicated algebra of
$A$.

\begin{lem}\label{lem 4.1} Let $T$ be an exceptional $\Am$-module having all
  projective-injec\-tive indecomposable modules as direct summands and
  let $M$ be an $\Am$-module. Assume that $M$ has a projective-injective
  injective envelope. Then a minimal left $\add T$-approximation of
  $M$ is a monomorphism.
\end{lem}  

\begin{pf}
Let $f_0:M\to T_0$ be a minimal left $\add T$-approximation and
$g:M\to I$ be an injective envelope. Since $I$ is
projective-injective, it lies in $\add T$, hence there exists
$h:T_0\to I$ such that $g=h\,f_0$. Since $g$ is injective, so is
$f_0$. 
\qed
\end{pf}  

\begin{cor}\label{cor 4.2}
An exceptional $\Am$-module $T$ is faithful if and only if it has all
projective-injective indecomposable $\Am$-modules as direct summands.
\end{cor}  

\begin{pf}
We have already shown the necessity. Conversely, assume  any
pro\-jec\-tive-injective indecomposable $\Am$-module to be a summand of
$T$. By Lemma \ref{lem 4.1}, a minimal left $\add T$-approximation of
$A_\Am$ is a mono\-mor\-phism. Therefore, there exists a mono\-mor\-phism
$\Am\!_\Am \monoto T_0$ with $T_0\in \add T$.
\qed
\end{pf}  
\end{subsection}

\begin{subsection}{Tilting modules}
We shall prove that, if $T$ is a faithful exceptional $\Am$-module, all
of whose non projective-injective summands lie in $\add \Lm$, then
there exists an $\Am$-module $X$ such that $T\oplus X $ is a tilting
module.

Clearly, if $\Am$ (or, equivalently, $A$) is representation-finite, then
this follows from \cite{RS}. We may then assume without loss of
generality that $A$ is representation-infinite.

\begin{lem}\label{lem 4.3}
Assume that $A$ is representation-infinite, that $T$ is  a faithful exceptional
$\Am$-module and that $M\in\Lm$. Then 
\begin{enumerate}
\item the injective envelope of $M$ is projective-injective, and
\item a minimal left $\add T$-approximation of $M$ is injective.
\end{enumerate}   
\end{lem}  

\begin{pf}
By Lemma \ref{lem 4.1}, it suffices to prove (1). Since, by Corollary
\ref{cor 3.5}, we have $M\in\Lm$ if and only if $M\le \zS_m$, it
suffices, by Lemma \ref{lem 3.1}, to show that $\zS_m$ contains no
$A_m$-module.
 We first assume that $m=2l$ is even, then 
\begin{equation}\nonumber
\zS_m\,=\, \zO_\Am^{-m}\,\zS_{0}\,=\,
\zO_\Am^{-2l}\,\zS_{0}\,=\,
\tau_\Am^{-l}\,(\nu_\Am^l\,\zS_{0}).
\end{equation}  
Since $\nu_\Am^{l}\,\zS_0$ is the set of projective $A_l$-modules, the
set $\zS_{m}$ consists of postprojective
$A_l$-modules (here, we are using the assumption that $\Am$ is
representation-infinite).  Since $l\ne
m$, this shows that 
 $\zS_m \cap \ind A_m = \emptyset$.
Assume now that $m=2l+1$ is odd, then 
\begin{equation}\nonumber
\zS_m\,=\, \zO_\Am^{-1}\,\zS_{2l}\,=\,
\zO_\Am^{-1}\,\tau_\Am^{-l}\,\nu_\Am^l\,\zS_{0}\,=\,
\tau_\Am^{-l}\,\zO_\Am^{-1} \,(\nu_\Am^l\,\zS_{0}).
\end{equation}  
But $\nu_\Am^{l}\,\zS_0$ is the set of all indecomposable
projective $A_l$-modules and, since $A$ is representation-infinite,
this  implies $\zO_\Am^{-1} 
\,(\nu_\Am^l\,\zS_{0})\subset \ind A_{l,l+1}$ (in the notation of
section \ref{section 2.1} and section \ref{section 3.1}). We thus get
again $\zS_m \cap \ind A_m = \emptyset$.
\qed
\end{pf}

\begin{rem} The statement of Lemma \ref{lem 4.3} is false in the
    representation-finite case: 
  let $A$ be the path algebra of the quiver $a\ot b$, then the
  Auslander-Reiten quiver of $A^{(1)}$ is:\\
\[
\xymatrix@!0@C=35pt@R=35pt{a_0\ar[dr]&&b_0\ar[dr]&&a_1\ar[dr]&&b_1\\
&{\begin{array}{c} \vspace{-10pt}b_0\\  a_0 \end{array}}\ar[dr]\ar[ur]
&&{\begin{array}{c} \vspace{-10pt}a_1\\b_0 \end{array}}\ar[r]\ar[ur]
&{\begin{array}{c} \vspace{-10pt}b_1\\\vspace{-10pt}a_1\\b_0 \end{array}}\ar[r]
&{\begin{array}{c} \vspace{-10pt}b_1\\a_1 \end{array}}\ar[ur]\\
&&{\begin{array}{c}\vspace{-10pt} a_1\\ \vspace{-10pt}b_0\\a_0 \end{array}}\ar[ur]
}
\]
Clearly, the simple module $a_1$ belongs to $\zS_1$ but its injective
envelope ${\begin{array}{c} \vspace{-12pt}b_1\\a_1 \end{array}} $
 is not projective.
The above proof fails here because the module 
${\begin{array}{c} \vspace{-12pt}b_0\\a_0 \end{array}}$
 is a projective-injective $A$-module and
  $\zO_\Am^{-1}{\begin{array}{c} \vspace{-12pt}b_0\\a_0 \end{array}}
 = a_1$ is not in $\ind A_{0,1}$. 
\end{rem}

\begin{prop}\label{prop 4.4}
Let $A$ be representation-infinite and let $T$ be a faithful
exceptional $\Am$-module such that $\pd T \le m$. Then there exists an
exact sequence 
\begin{equation}\nonumber
0\longrightarrow A \stackrel{f_0}{\longrightarrow} T_0  \stackrel{f_1}{\longrightarrow} T_1
\stackrel{f_2}{\longrightarrow} \cdots  \stackrel{f_{m-1}}{\longrightarrow} T_{m-1}
\stackrel{f_m}{\longrightarrow} L_m    \longrightarrow 0  
\end{equation}
in $\textup{mod}\, \Am$ such that 
\begin{enumerate}
\item $T_i\in\add T$ for all $i$,
\item $L_i=\Coker f_{i-1}$ lies in $\add\mathcal{L}_i(\Am)$, and
\item each of the induced monomorphisms $L_i\monoto T_i$ is a minimal
  left $\add T$-approximation.
\end{enumerate}  
\end{prop}  

\begin{pf}
We construct this sequence by induction on $s<m$. First, since $T$ is faithful
and $\pd A_\Am=0$, it follows from Lemma \ref{lem 4.3} that a minimal
left $\add T$-approximation of $A_\Am$ is a monomorphism, and thus we
have a short exact sequence 
\begin{equation}\nonumber
0\longrightarrow A\stackrel{f_0}{\longrightarrow} T_0 \stackrel{g_1}\longrightarrow L_1 \longrightarrow 0
\end{equation}  
in $\textup{mod}\, \Am$. 

Assume now that $s<m$ and that we have an exact
sequence 
\begin{equation}\nonumber
0\longrightarrow A\stackrel{f_0}{\longrightarrow} T_0 \stackrel{f_1}{\longrightarrow} T_1 \longrightarrow \cdots
\stackrel{f_s}{\longrightarrow} T_s 
\end{equation}  
such that $T_i \in\add T, \ L_i=\Coker f_{i-1}\in\add\mathcal{L}_i(\Am)$
and each of the 
induced monomorphisms $f_i':L_i\monoto T_i$ is a minimal
  left $\add T$-approximation.
Let $L_{s+1}=\Coker f_s$. 
If $L_{s+1}=0$, then our sequence stops and
it is of the form
\begin{equation}\nonumber
0\longrightarrow A {\longrightarrow} T_0 \longrightarrow T_1 \longrightarrow \cdots
{\longrightarrow} T_s \longrightarrow 0
\end{equation}  
and all $L_i=0$, $T_i=0$ if $s<i\le m$. This sequence clearly
satisfies conditions (1)-(3).
If $L_{s+1}\ne 0$, consider the short exact sequence 
\begin{equation}\nonumber
0\longrightarrow L_s\stackrel{f'_s}{\longrightarrow} T_s
\stackrel{g_{s+1}} \longrightarrow L_{s+1} \longrightarrow 0
\end{equation}  
in $\textup{mod}\, \Am$. By the induction hypothesis, $L_s\in
\mathcal{L}_s(\Am)\subset\Lm  $,
therefore, by Lemma \ref{lem 4.3}, $L_s $ has a projective-injective
injective envelope $I$, so that we have a commutative diagram with exact rows
\[
\xymatrix{
0\ar[r]&L_s \ar[r]^{f'_s} \ar[d]_{\textup{id}}
&T_s \ar[r]^{g_{s+1}} \ar[d]
&L_{s+1} \ar[r] \ar[d]&0\\
0\ar[r]&L_s \ar[r]
&I \ar[r]
&\zO_\Am^{-1}\,L_{s} \ar[r] &0 .
}
\]
Since $L_s\in\mathcal{L}_s(\Am)$ implies $L_s\le \zS_s$, we have
\begin{equation}\nonumber
L_{s+1}\le\zO_\Am^{-1}\,L_{s}\le
  \zO_\Am^{-1}\,\zS_{s} =\zS_{s+1}\le\zS_m 
\end{equation}  
(because $s+1\le m$). In particular $L_{s+1}\in
\mathcal{L}_{s+1}(\Am)$. By Lemma \ref{lem 4.3}, a minimal left 
$\add T$-approximation $f'_{s+1}:L_{s+1}\to T_{s+1}$ is a monomorphism
and in particular $T_{s+1}\ne 0$.
We obtain the exact sequence
\begin{equation}\nonumber
0\longrightarrow A \stackrel{f_0}{\longrightarrow} T_0  \stackrel{f_1}{\longrightarrow} T_1
{\longrightarrow} \cdots  \stackrel{f_{s}}{\longrightarrow} T_{s}
\stackrel{f_{s+1}}{\longrightarrow} T_{s+1} \longrightarrow
L_{s+1}\longrightarrow 0
\end{equation}
where $f_{s+1}= f_{s+1}'\circ g_{s+1}$,
 and this completes the proof.
\qed
\end{pf}  
 
\begin{thm}\label{thm 4.5} Let $A$ be a hereditary algebra, $m\ge 1$
  and $\Am$ be the $m$-th replicated algebra of $A$. If $T$ is a faithful
  exceptional $\Am$-module with $\pd T \le m$, then there exists an
  $\Am$-module $X$ such that $T\oplus X$ is a tilting $\Am$-module and
  $\pd X \le m$.
\end{thm}  

\begin{pf}
We may, by \cite{RS}, assume that $\Am$, or equivalently, $A$, is
repre\-sen\-tation-infinite. Then, by Proposition \ref{prop 4.4}, there
exists an exact sequence
\begin{equation}\nonumber
0\longrightarrow A \stackrel{f_0}{\longrightarrow} T_0  \stackrel{f_1}{\longrightarrow} T_1
{\longrightarrow} \cdots  \stackrel{f_{m-1}}{\longrightarrow} T_{m-1}
\stackrel{f_m}{\longrightarrow} L_m   \longrightarrow 0 
\end{equation}
in $\textup{mod}\, \Am$ such that 
\begin{enumerate}
\item $T_i\in\add T$ for all $i$,
\item $L_i=\Coker f_{i-1}$ lies in $\add\mathcal{L}_i(\Am)$ (or,
 equivalently, $\pd L_i\le i$), and
\item each of the induced monomorphisms $L_i\monoto T_i$ is a minimal
  left $\add T$-approximation.
\end{enumerate}  
Actually, we have one of two cases. If there exists $p\le m$ such that
$L_p=0$, then the above sequence reduces to
\begin{equation}\nonumber
0\longrightarrow A \stackrel{f_0}{\longrightarrow} T_0  \stackrel{f_1}{\longrightarrow} T_1
{\longrightarrow} \cdots  \stackrel{f_{p}}{\longrightarrow} T_{p} \longrightarrow 0
\end{equation}
(with $T_i\in\add T$ for all $i$) and then $T$ is clearly a tilting
module.
If not we may assume $L_m\ne 0$. 
We then prove by induction on $s$, with $1\le s\le m$, that for all $i\ge 1$
\begin{eqnarray}
&& \textup{$\Ext_\Am^i(L_s,T)=0$.}\label{a}\\
 &&\textup{$\Ext_\Am^i(T,L_s) \isomorphe   \Ext_\Am^{i+s}(T,A) $.}\label{b}\\
&&\textup{$\Ext_\Am^i(L_s,L_s) \isomorphe   \Ext_\Am^{i+s}(L_s,A) $.}\label{c}
\end{eqnarray}     
Assume first that $s=1$, and apply $\Hom_\Am(-,T)$ to the short exact
sequence 
\begin{equation}\label{stern}
0\longrightarrow A \stackrel{f_0}{\longrightarrow} T_0 \longrightarrow L_1 \longrightarrow 0.
\end{equation}  
This yields an exact sequence 
\[\begin{array}{rcccl}
\Hom_\Am(T_0,T)&\stackrel{\Hom\!_\Am(f_0,T)}{\longrightarrow} &
\Hom_\Am(A,T)&\longrightarrow&\Ext_\Am^1(L_1,T)\\
&\longrightarrow&\Ext_\Am^1(T_0,T)&=&0
\end{array}\]  
where the last equality follows from the exceptionality of $T$.
Since $f_0$ is a left $\add T$-approximation, $\Hom_\Am(f_0,T)$ is
surjective. Hence $\Ext_\Am^1(L_1,T)=0$. This implies (\ref{a}), because $\pd
L_1\le 1$.

Applying now $\Hom_\Am(T,-)$ to (\ref{stern}) yields the exact sequence
\begin{equation}\nonumber
0=\Ext_\Am^i(T,T_0)\to\Ext_\Am^i(T,L_1)\to
\Ext_\Am^{i+1}(T,A)\to\Ext_\Am^{i+1}(T,T_0)=0,
\end{equation} 
hence $\Ext_\Am^i(T,L_1)\isomorphe\Ext_\Am^{i+1}(T,A) $. This is just (\ref{b}).

Finally, applying $\Hom_\Am(L_1,-)$ to (\ref{stern}) yields the exact sequence
\[\begin{array}{rcccl}
0=\Ext_\Am^1(L_1,T_0)&\longrightarrow&\Ext_\Am^1(L_1,L_1)&\longrightarrow&
\Ext_\Am^{2}(L_1,A)\\
&\longrightarrow&\Ext_\Am^{2}(L_1,T_0)&=&0 
\end{array} \]
hence $\Ext_\Am^1(L_1,L_1)\isomorphe\Ext_\Am^{2}(L_1,A) \,(=0) $, and
thus (\ref{c}), because $\pd L_1\le 1$.

Let now $s>1$ and consider the short exact sequence 
\begin{equation}\label{stern2}
0\longrightarrow L_{s-1} \stackrel{f'_{s-1}}{\longrightarrow} T_{s-1} \longrightarrow L_s \longrightarrow 0.
\end{equation} 
Applying $\Hom_\Am(-,T)$ yields an exact sequence 
\[
\begin{array}{ccccl}
\Hom_\Am(T_{s-1},T)&\stackrel{\Hom_\Am(f'_{s-1},T)}{\longrightarrow}& 
\Hom_\Am(L_{s-1},T)&\longrightarrow&\Ext_\Am^1(L_s,T)\\
&\longrightarrow&
\Ext_\Am^1(T_{s-1},T)&=&0
\end{array}  \]
Since $f'_{s-1}$ is a left $\add T$-approximation, we get
$\Ext_\Am^1(L_s,T)=0$. 
The same long exact sequence yields, for $i\ge 2$,
\begin{equation}\nonumber
\Ext_\Am^i(L_s,T)\isomorphe \Ext_\Am^{i-1}(L_{s-1},T)=0
\end{equation} 
where the second equality follows from the induction hypothesis. This
shows (\ref{a}).

Applying  $\Hom_\Am(T,-)$ to (\ref{stern2}) yields the exact sequences
\[\begin{array}{ccccl}
0=\Ext_\Am^i(T,T_{s-1})&\longrightarrow&\Ext_\Am^i(T,L_s)&\longrightarrow&
\Ext_\Am^{i+1}(T,L_{s-1})\\
&\longrightarrow& \Ext_\Am^{i+1}(T,T_{s-1})&=&0  
\end{array} \]
thus $\Ext_\Am^i(T,L_{s})\isomorphe\Ext_\Am^{i+1}(T,L_{s-1})$ and
$\Ext_\Am^{i+1}(T,L_{s-1})\isomorphe\Ext_\Am^{i+s}(T,A) $, by the
induction hypothesis. This 
shows  (\ref{b}).

Finally, applying $\Hom_\Am(L_s,-)$ to (\ref{stern2}) yields the exact
sequences 
\[\begin{array}{ccccl}
0=\Ext_\Am^i(L_s,T_{s-1})&\longrightarrow&\Ext_\Am^i(L_s,L_s)&\longrightarrow&
\Ext_\Am^{i+1}(L_s,L_{s-1})\\
&\longrightarrow&\Ext_\Am^{i+1}(L_s,T_{s-1})&=&0 
\end{array}\]
hence $\Ext\!_\Am^i(L_s,L_s)\isomorphe\Ext\!_\Am^{i+1}(L_s,L_{s-1})$.
Similarly, applying $\Hom\!_\Am(L_s,-)$ to the short exact sequence
\begin{equation*}
0\longrightarrow L_{s-2} \stackrel{f'_{s-2}}{\longrightarrow} T_{s-2}
\longrightarrow L_{s-1} \longrightarrow 0. 
\end{equation*} 
yields $\Ext_\Am^{i+1}(L_s,L_{s-1})\isomorphe\Ext_\Am^{i+2}(L_s,L_{s-2})$.
Continuing in this way, one gets eventually
$\Ext_\Am^i(L_s,L_s)\isomorphe \Ext_\Am^{i+s}(L_s,L_0)$
with $L_0=A$. This
shows  (\ref{c})  and
completes the proof of  our claim.

Let now $s=m$. We deduce that, for all $i\ge 1$ 
\[
\begin{array}{ll}
\Ext_\Am^i(L_m,T)=0. \\
\Ext_\Am^i(T,L_m) \isomorphe   \Ext_\Am^{m+i}(T,A) =0  &\textup{since
  $\pd T\le m$. }\\
\Ext_\Am^i(L_m,L_m) \isomorphe   \Ext_\Am^{m+i}(L_m,A)=0 
&\textup{since $\pd L_m\le m$.}\\
\end{array}
\] 
This shows that $T\oplus L_m$ is a tilting $\Am$-module, and completes
the proof of the theorem.
\qed
\end{pf}

\begin{rem} Observe that $T$ has usually many possible complements. 
 Our proof constructs only one.
\end{rem}  

\begin{cor}\label{cor 4.6}
Let $A$ be a hereditary algebra, $m\ge 1$ and $\Am$ be the
$m$-replicated algebra of $A$. Let $T$ be a basic, faithful,
exceptional $\Am$-module with $\pd T\le m$ and such that the number of
indecomposable summands of $T$ is equal to the rank of the
Grothendieck group of $\Am$.
Then $T$ is a tilting $\Am$-module.
\end{cor}  
\end{subsection}
\end{section}

\begin{section}{Tilting modules and tilting objects}

\begin{subsection}{Main result}
Let $A$ be a hereditary algebra, $m\ge 1$ and $X$ be an object in the
$m$-cluster category $\CmA$.  Then  $X$ is said to be {\em basic} (or {\em
  multiplicity free}) if, whenever $X=\oplus_{i=1}^t X_i$ where all the
  $X_i$ are indecomposable, we have $X_i\ne X_j$ for $i\ne j$.
The object 
$X$ is called {\em exceptional} if $\Ext^i_{\CmA}(X,X)=0$ for all $i$
with $1\le i\le m$ and it is called {\em tilting} if it is exceptional
and the   number of isomorphism classes of its indecomposable summands
is equal to the rank of the 
  Grothendieck group of $A$ (see \cite{T}).

Let now $T$ be an exceptional $\Am$-module.
 Then we can always write $T$ in the form $T=T'\oplus P$, where
$P$ is projective-injective and $T'$ has no projective-injective
indecomposable summands. We say that $T$ is an {\em $\mathcal{L}_m$-exceptional
  module} if $T'\in\add\Lm$ (or, equivalently, if $\pd T \le m$).

In this section, we always assume our exceptional objects and modules
to be basic. By abuse of notation, modules will often be denoted by
the same letter even when considered as objects in different
categories. We now prove Theorem \ref{thmC}.

  \begin{thm}\label{thm 5.1}
   Let $A$ be a hereditary algebra    and $\Am$  be its $m$-th
   replicated  algebra. There is a one-to-one
   correspondence: 
\begin{equation}\nonumber
 \{\text{basic } \mathcal{L}_{m}-\text{exceptional modules}\}\leftrightarrow
 \{\text{basic exceptional objects in }\mathcal{C}_{m} (A) \},
\end{equation} 
 which is given by
 $T=T'\oplus P \mapsto \pi (T')$. 
 \end{thm}
 \begin{pf}
Let $T=T'\oplus P$ be a basic $ \mathcal{L}_{m}$-{exceptional
 module}. We claim that
 $X=\pi(T')$ is an exceptional object in $\mathcal{C}_{m} (A)$,
 that is,
 $\text{Ext}^i_{\CmA}(X,X) = 0$ for  $1 \le i \le m.$
 By the definition of $\text{Ext}^i_{\CmA}$ this amounts to proving that
  \begin{equation} \label{homzero}
  \text{Hom}_{D^b(\text{mod} A)}(X_x,\tau^{-s}X_y[ms+i]) = 0
  \end{equation}
(where $\tau=\tau_{\DA}$)
   for all $s \in \Zset$, all $i$ such that $1 \le i \le m$ and all
    indecomposable summands $X_x,X_y$ of $X$. Denote by $T_x$ and
    $T_y$ the indecomposable $\Am$-modules lying in $\Lm$ which
    correspond to $X_x$ and $X_y$, respectively.
 We show equation (\ref{homzero}) by distinguishing various cases
    according to the value of $s$.
\begin{enumerate}
\item  If $s=0$ then the equation (\ref{homzero}) holds for all $i$
  since $T$ is an $\mathcal{L}_{m}-$exceptional module. 
\item If $s=-1$ then we have for all $i$
\begin{eqnarray}
\text{Hom}_{D^b(\text{mod} A)}(X_x,\tau X_y[-m+i])  &\isomorphe& 
\sHom_{\hat{A}}(T_x,\tau_{\hat{A}}\, \Omega^{m-i}_{\hat{A}}\, T_y)\label{11}\\
  &\isomorphe& D\Ext^1_{\hat{A}}(\Omega^{m-i}_{\hat{A}}\, T_y, T_x)\label{12} \\
  &\isomorphe& D\text{Ext}^{m+1-i}_{\hat{A}} (T_y, T_x)\label{13}\\
 &\isomorphe& D\text{Ext}^{m+1-i}_{\Am} (T_y, T_x)\label{14}\\
&=&0 \label{15}
\end{eqnarray}  
where (\ref{11}) follows from Theorem \ref{Happelthm},
(\ref{12}) is the Auslander-Reiten formula in $\DA$, (\ref{13}) follows from the
definition of $\zO$, (\ref{14}) holds because $T\in \add\Lm$ and
(\ref{15})  because $T$ is  exceptional.
\item If $s \le -2$ then $ms+i \le -2m+i \le -m$, so $X_y[ms+i]$ lies
  in some  $(\ind A)[j] $ with $j<0$ except in the case where $s=-2,\ i=m$
  and $X_y \in (\ind A)[m]$, in which case  $X_y[ms+i] \in (\ind A)[0]$,
  and then
  $\tau^{-s}X_y[ms+i](=\tau^2 X_y[ms+i])$ lies in some  $(\ind A)[j] $ with $j<0$.
 In either case, there are no non-zero morphisms from
  $X_x$ to $\tau^{-s}X_y[ms+i]$. 
\item 
 If $s \ge 1$, then $ms+i\ge m+i\ge m$ so the only possibility to have
 a non-zero morphism from  
 $X_x$ to $\tau^{-s}X_y[ms+i]$ is when $s=1,i=1, X_x =P[m]$ for some
 indecomposable projective $A$-module $P$, 
 and $\tau^{-1}\, X_y \in (\ind A)[0]$. But then 
 $\tau^{-1}\, X_y[ms+i] \in (\ind A)[m+1]$ and there is no non-zero
 morphism from $P[m]$ to $(\ind A)[m+1]$.
\end{enumerate}   
 
  Conversely, assume that $X=\pi(T')$ is a basic exceptional
 object in $\CmA$. We claim that $T'$ is a basic
 $\mathcal{L}_{m}-$exceptional $\Am$-module.
Clearly,  $T'$ is basic and of finite projective dimension. Moreover $T' \in
 \mathcal{L}_{m}(\Am)$, by Theorem \ref{thm2}. 
 Suppose that there exist indecomposable summands $T_x,T_y$ of $T'$
 such that $\text{Ext}^i_{\Am}(T_x,T_y) \neq 0$ for some $i$ with $1\le
 i\le m$. Then 
\[ \begin{array}{rcccl}
\sHom _{\hat{A}}(T_y,\tau_{\Ahat}\Omega_{\Ahat}^{i-1}\,T_x)&\isomorphe&
D\text{Ext}^1_{\hat{A}}(\Omega^{i-1}_{\hat{A}}T\,_x,T_y)&\isomorphe&
D  \text{Ext}^i_{\hat{A}}(T_x,T_y)\\
&\isomorphe & D \text{Ext}^i_{\Am}(T_x,T_y)
\end{array}\]
implies that   
\begin{equation}\nonumber
\Hom_{\DA}(X_y,\tau X_x[-i+1]) \neq 0
\end{equation}  
(where $X_x,\,X_y$ denote, as before, the indecomposable summands of
$X$ which correspond to $T_x,\,T_y$, respectively). Thus we have
\begin{eqnarray*}
 \text{Ext}^i_{\CmA}(X_x,X_y)&\isomorphe&
 \text{Ext}^1_{\CmA}(X_x[-i+1],X_y)\\
&\isomorphe& \text{Hom} _{\CmA}(X_y,\tau_{\CmA}X_x[-i+1])\\
&\ne& 0
\end{eqnarray*}  
contradicting the hypothesis that $X$ is an exceptional object in
$\CmA$
\qed 
 \end{pf}

\begin{cor} Let $A$ be a hereditary algebra   and $\Am$  be its $m$-th
   replicated  algebra. There is a one-to-one
  correspondence
\begin{equation}\nonumber
 \{\text{basic } \mathcal{L}_{m}-\text{tilting modules}\}\leftrightarrow 
 \{\text{basic tilting objects in }\mathcal{C}_{m} A \},
\end{equation} 
 which is given by
 $T=T'\oplus P \mapsto \pi (T')$. 
\end{cor}  

\begin{pf}
Assume $T$ is a basic $\mathcal{L}_{m}$-tilting module,
then $T=T'\oplus P$ where $P$ has $nm$ indecomposable summands (here
$n$ is the rank of the Grothendieck group of $A$) while $T'$ has $n$
indecomposable summands. But then $X=\pi (T')$ has also $n$
indecomposable summands. Since, by Theorem \ref{thm 5.1}, $X$ is
exceptional, it is tilting.

Conversely, if $X$ is a tilting object in $\CmA$ then it has $n$
indecomposable summands. Let $T=T'\oplus P$ where $P$ is the direct
sum of all projective-injective indecomposable $\Am$-modules and $T'$ is
such that $\pi(T')=X$. Since $P$ has $nm$ indecomposable
summands, and $T'$ has $n$, then $T$ has $nm+n$ indecomposable
summands. But $nm+n$ is equal to the rank of the Grothendieck group of
$\Am$ and, by Corollary \ref{cor 4.6}, $T$ is a tilting module.
\qed
\end{pf}  
\end{subsection}

\begin{subsection}{Application to the case $m=1$}
In \cite{ABST1}, it is shown that there is a bijection between the
tilting objects in the cluster category $\CA=\mathcal{C}_1(A)$ and the
$\mathcal{L}$-tilting modules in the duplicated algebra
$\Abar=A^{(1)}$. Now, recall that in this context, tilting modules
are of projective dimension $1$. By Corollary \ref{cor 3.5}, this is
equivalent to the fact that the non projective-injective
indecomposable summands lie in
$\mathcal{L}_{\Abar}=\mathcal{L}_1({\Abar})$. Therefore, we have the
following corollary.
\begin{cor}
Let $A$ be a hereditary algebra. There is a one-to-one
  correspondence
\begin{equation}\nonumber
 \{\text{basic tilting $\Abar$-modules}\}\leftrightarrow
 \{\text{basic cluster tilting objects in }\mathcal{C}(A) \},
\end{equation} 
 which is given by
 $T=T'\oplus P \mapsto \pi (T')$. 
\end{cor}  
\end{subsection}

\begin{subsection}{Example} 
  Let $A$ be given by the quiver 
\[\xymatrix{1&2\ar[l]&3\ar[l]&4\ar[l]} \]
Then the Auslander-Reiten quiver of $A^{(2)}$ is given by
\[
\xy
(-14,0)*++{\xy(4,0)*{}="A"; (-40,0)*{}="B"; (24,40)*{}="C";(-20,40)*{}="D";
"A"; "B" **\dir{-};
"D"; "C" **\dir{-};
"C"; "A" **\dir{-}; 
"D"; "B" **\dir{-}; 
\endxy};
(-20,9)*++{
\xymatrix@R=10pt@C=0pt
{
&&&&&&\circ\ar[rdd]&&&&\circ\ar[rdd]&&&&\circ\ar[rdd]&&&&\circ\ar[rdd]&&&&\\
&&&&\circ\ar[rd]&&&&\circ\ar[rd]&&&&\circ\ar[rd]&&&&\circ\ar[rd]&&&&&\\
&&&\diamond\ar[ru]\ar[rd]&&\bullet\ar[ruu]\ar[rd]&&
\bullet\ar[ru]\ar[rd]&&\diamond\ar[ruu]\ar[rd]&&
\bullet\ar[ru]\ar[rd]&&\bullet\ar[ruu]\ar[rd]
&&\bullet\ar[ru]\ar[rd]&&\bullet\ar[ruu]\ar[rd]&&\bullet\ar[rd]\\
&&\bullet\ar[ru]\ar[rd]&&\bullet\ar[ru]\ar[rd]&&
\bullet\ar[ru]\ar[rd]&&\bullet\ar[ru]\ar[rd]&&
\bullet\ar[ru]\ar[rd]&&\bullet\ar[ru]\ar[rd]&&\bullet\ar[ru]\ar[rd]
&&\bullet\ar[ru]\ar[rd]&&\bullet\ar[ru]\ar[rd]&&\bullet\ar[rd]\\
&\bullet\ar[ru]\ar[rd]&&\bullet\ar[ru]\ar[rd]&&
\bullet\ar[ru]\ar[rd]&&
\bullet\ar[ru]\ar[rd]&&\bullet\ar[ru]\ar[rd]&&
\bullet\ar[ru]\ar[rd]&&\bullet\ar[ru]\ar[rd]&&\bullet\ar[ru]\ar[rd]
&&\bullet\ar[ru]\ar[rd]&&\bullet\ar[ru]\ar[rd]&&\bullet\ar[rd]\\
\diamond\ar[ru]&&\bullet\ar[ru]&&
\bullet\ar[ru]&&\diamond\ar[ru]&&
\bullet\ar[ru]&&\bullet\ar[ru]&&
\bullet\ar[ru]&&\bullet\ar[ru]&&\bullet\ar[ru]
&&\bullet\ar[ru]&&\bullet\ar[ru]&&\bullet\\
}};
\endxy
\]

where we have indicated the 2-left part $\mathcal{L}_2 (A^{(2)})$.  We
have also indicated an 
$\mathcal{L}_2$-tilting module $T=T'\oplus P$, where $T' \in
\add \mathcal{L}_2$. The summands of $T'$ are indicated by diamonds and the
(projective-injective) summands of $P$ by circles.

\end{subsection}
\end{section}

\end{document}